# Residue symbols and Jantzen-Seitz partitions

C. Bessenrodt and J.B. Olsson

Abstract. Jantzen-Seitz partitions are those $p$-regular partitions of $n$ which label $p$-modular irreducible representations of the symmetric group $S_n$ which remain irreducible when restricted to $S_{n-1}$; they have recently also been found to be important for certain exactly solvable models in statistical mechanics. In this article we study their combinatorial properties via a detailed analysis of their residue symbols; in particular the $p$-cores of Jantzen-Seitz partitions are determined.

## 1 Introduction

The Mullineux symbols are combinatorial objects which were introduced in order to understand the following question in the modular representation theory of the symmetric groups $S_n$, $n$ a natural number.

For a given prime $p$ the $p$-modular irreducible representations $D^\lambda$ of $S_n$ are labeled in a canonical way by the $p$-regular partitions $\lambda$ of $n$. When the modular irreducible representation $D^\lambda$ of $S_n$ is tensored by the sign representation we get a new modular irreducible representation $D^{\lambda^P}$. The question about the connection between the $p$-regular partitions $\lambda$ and $\lambda^P$ was answered in 1995 by the proof of the so-called "Mullineux Conjecture".

The conjecture was formulated in 1979, when Mullineux [15] defined a bijective involutory map $\lambda \to \lambda^M$ on the set of $p$-regular partitions and conjectured that this map coincides with the map $\lambda \to \lambda^P$. The statement "$M = P$" is the Mullineux conjecture.

The map $M$ is an operation on the Mullineux symbols of $p$-regular partitions as described below. The proof of the Mullineux conjecture was made possible by a series of papers by A. Kleshchev [9]-[11] on "modular branching", i.e. on the restrictions of modular irreducible representations from $S_n$ to $S_{n-1}$. Kleshchev's result reduced the proof of the Mullineux Conjecture to a purely combinatorial statement about the compatibility of the Mullineux map with the removal of "good nodes" (see below). A lengthy proof of this combinatorial statement was then given in a paper by Ford and Kleshchev [5]. By introducing a variation of the Mullineux symbols, called *residue symbols* the present authors were able to give a shorter and more transparent proof [2]. One advantage of the residue symbols is the easy recognition of good and normal nodes in the corresponding $p$-regular partitions. Thus they are useful tools when branching and Mullineux conjugation both have to be under control.

In this paper we study further properties of the residue symbols. The motivation for most of our results is the following. We want to study the analogues of Kleshchev's branching results for alternating groups, in particular the question of which modular irreducible representations of $A_n$ remain irreducible upon restriction to $A_{n-1}$ (our results on branching in alternating groups will appear in a forthcoming paper). In section 2 we state as motivation some facts about modular branching and Jantzen-Seitz partitions; these partitions label those modular representations of $S_n$ which remain irreducible upon restriction to $S_{n-1}$ and it is necessary to understand which Jantzen-Seitz partitions are fixed under the Mullineux map. We explain the relevance of residue symbols of Jantzen-Seitz partitions. The residue symbols of these partitions have a very simple structure (see section 3); they may be described as paths in a special directed graph. In particular, we use this description to show that the $p$-core of such a partition is a rectangular partition determined by the end of the path in the diagram; this problem has also come up in the context of studying special models in statistical mechanics [4]. Moreover, we describe explicitly the changes in the weights of the Jantzen-Seitz partitions constructed along the path. A similar analysis is applied in section 4 to Jantzen-Seitz partitions fixed under the Mullineux map; their $p$-cores are described and also the possible weights for such partitions with a given $p$-core are determined.

# 2 Preliminaries

We use two types of notation to describe a partition $\lambda$ of $n$. The first is the usual one:

$$\lambda = (\lambda_1 \geq \lambda_2 \geq \ldots \geq \lambda_k > 0)$$

where the parts $\lambda_i \in \mathbb{N}$ satisfy $\lambda_1 + \lambda_2 + \ldots + \lambda_k = n$. Alternatively we use an "exponential" notation

$$\lambda = (l_1^{a_1}, l_2^{a_2}, \ldots, l_t^{a_t})$$

where $l_i$, $a_i \in \mathbb{N}$ for all $i$, $l_1 > l_2 > \ldots > l_t$ and $\sum_{i=1}^t a_i l_i = n$. In this notation, the first $a_1$ parts of $\lambda$ equal $l_1$, the next $a_2$ parts of $\lambda$ equal $l_2$ etc.

**Example.** $\lambda = (5,5,4,1,1,1) = (5^2, 4, 1^3)$.

Let $p \geq 2$ be an integer. Whenever $p$ occurs in statements about representations, we assume $p$ to be a prime integer. In all other connections $p$ may be arbitrary.

A partition $\lambda = (l_1^{a_1}, l_2^{a_2}, \ldots, l_t^{a_t})$ is called *$p$-regular*, if $1 \leq a_i \leq p-1$ for $i = 1, \ldots, t$. We will consider only $p$-regular partitions in this paper.

Let $\lambda = (\lambda_1 \geq \lambda_2 \geq \ldots \geq \lambda_k > 0)$ be a partition of $n$. Then

$$Y(\lambda) = \{(i,j) \in \mathbb{Z} \times \mathbb{Z} \mid 1 \leq i \leq k, 1 \leq j \leq \lambda_i\} \subset \mathbb{Z} \times \mathbb{Z}$$

is the Young diagram of $\lambda$, and $(i,j) \in Y(\lambda)$ is called a *node* of $\lambda$. If $A = (i,j)$ is a node of $\lambda$ and $Y(\lambda) \setminus \{(i,j)\}$ is again a Young diagram of a partition, then $A$ is called a *removable* node and $\lambda \setminus A$ denotes the corresponding partition of $n-1$.

Similarly, if $A = (i,j) \in \mathbb{N} \times \mathbb{N}$ is such that $Y(\lambda) \cup \{(i,j)\}$ is the Young diagram of a partition of $n+1$, then $A$ is called an *indent* node of $\lambda$.

The *p-residue* of a node $A = (i,j)$ is defined to be the residue modulo $p$ of $j-i$, denoted res $A = j - i \pmod p$. The $p$-residue diagram of $\lambda$ is obtained by writing the $p$-residue of each node of the Young diagram of $\lambda$ in the corresponding place.

**Example.** Take $p = 5$, $\lambda = (6^2, 5, 4)$. Then the 5-residue diagram of $\lambda$ is

$$\begin{array}{cccccc} 0 & 1 & 2 & 3 & 4 & 0 \\ 4 & 0 & 1 & 2 & 3 & 4 \\ 3 & 4 & 0 & 1 & 2 & \\ 2 & 3 & 4 & 0 & & \end{array}$$

We write our partition $\lambda$ in the form $\lambda = (l_1^{a_1}, l_2^{a_2}, \ldots, l_t^{a_t})$ and define for $1 \leq i \leq j \leq t$:

$$\beta(i,j) = l_i - l_j + \sum_{k=i}^{j} a_k \quad \text{and} \quad \gamma(i,j) = l_i - l_j + \sum_{k=i+1}^{j} a_k$$

Moreover, for $i \in \{1, \ldots, t\}$ let

$$M_i = \{j \mid 1 \leq j < i, \beta(j,i) \equiv 0 \pmod p\}$$

We then call the $i$'th corner node *normal* if and only if for all $j \in M_i$ there exists $d(j) \in \{j+1, \ldots, i-1\}$ satisfying $\beta(j, d(j)) \equiv 0 \pmod p$, and such that $|\{d(j) \mid j \in M_i\}| = |M_i|$.

We call $i$ *good* if it is normal and if $\gamma(i,i') \not\equiv 0 \pmod p$ for all normal $i' > i$.

If $i$ is normal (resp. good), we call the removable node in the $i$th block of $\lambda$ *normal* (resp. *good*). These properties of removable nodes may be easily read off the $p$-residue diagram of $\lambda$. One sees immediately that $\beta(i,j)$ is just the length of the path from the node at the beginning of the $i$-th block of $\lambda$ to the node at the end of the $j$th block of $\lambda$. The condition $\beta(i,j) \equiv 0 \pmod p$ is then equivalent to the equality of the $p$-residue of the indent node in the outer corner of the $i$th block and the $p$-residue of the removable node at the inner corner of the $j$th block.

Similarly, $\gamma(i,j) \equiv 0 \pmod p$ is equivalent to the equality of the $p$-residues of the removable nodes at the end of the $i$th and $j$th block.

We will say that a node $A = (i,j)$ is *above* the node $B = (i',j')$ (resp. $B$ is *below* $A$) if $i < i'$, and write this relation as $B \nearrow A$. Then a removable node $A$ of $\lambda$ is *normal* if for any $B \in \mathcal{M}_A = \{C \mid C \text{ indent node of } \lambda \text{ above } A \text{ with res } C = \text{res } A\}$ we can choose a removable node $C_B$ of $\lambda$ with $A \nearrow C_B \nearrow B$ and res $C_B =$ res $A$, such that $|\{C_B \mid B \in \mathcal{M}_A\}| = |\mathcal{M}_A|$. A node $A$ is *good* if it is the lowest normal node of its $p$-residue.

The *height* of a normal node $A$ is defined as

$$\mathrm{ht}(A) = 1 + |\{B \mid B \text{ normal node of } \lambda \text{ with } A \nearrow B \text{ and res } A = \text{res } B\}|$$

Let us look at the example $\lambda = (12, 7^2, 5^3, 3, 1^3)$.
In the 5-residue diagram below we have included also the indent nodes (which do not belong to $Y(\lambda)$). The equality of the residues of the indent node at the beginning of the second block and of the removable node in the fourth block corresponds to the fact that $\beta(2,4) \equiv 0 \pmod 5$. From the diagram it is also obvious that $\gamma(4,5) \equiv 0 \pmod 5$. In the diagram we have marked the normal nodes in boldface and noted their height.

| | | | | | | | | | | | | | |
|---|---|---|---|---|---|---|---|---|---|---|---|---|---|
| 0 | 1 | 2 | 3 | 4 | 0 | 1 | 2 | 3 | 4 | 0 | **1** | 2 | height 1 |
| 4 | 0 | 1 | 2 | 3 | 4 | 0 | 1 | | | | | | |
| 3 | 4 | 0 | 1 | 2 | 3 | **4** | | | | | | | height 1 |
| 2 | 3 | 4 | 0 | 1 | 2 | | | | | | | | |
| 1 | 2 | 3 | 4 | 0 | | | | | | | | | |
| 0 | 1 | 2 | 3 | **4** | | | | | | | | | height 2 |
| 4 | 0 | 1 | 2 | | | | | | | | | | |
| 3 | 4 | | | | | | | | | | | | |
| 2 | | | | | | | | | | | | | |
| **1** | | | | | | | | | | | | | height 2 |
| 0 | | | | | | | | | | | | | |

The set $M_i$ corresponds in this picture to taking the removable node, say $A$, at the end of the $i$th block and then collecting into $M_i$ (resp. $\mathcal{M}_A$) all the indent nodes above this block of the same $p$-residue as $A$. For $i$ resp. $A$ being normal, we then have to check whether for any such indent node, $B$ say, at the end of the $j$th block we can find a removable node $C = C_B$ between $A$ and $B$ of the same $p$-residue, and such that the collection of all these removable nodes has the same size as $M_i$ (resp. $\mathcal{M}_A$). The node $A$ (resp. $i$) is then good if $A$ is the lowest normal node of its $p$-residue.

The critical condition for the normality of $i$ resp. $A$ above is just a lattice condition: it says that in any section above $A$ there are at least as many removable nodes of the $p$-residue of $A$ as there are indent nodes of the same residue.

A convenient way to use this lattice condition to find normal nodes and their height is to use the node sequence $N(\lambda)$ of $\lambda$ as defined below. First we introduce signature sequences.

A $(p)$-*signature* is a pair $c\varepsilon$ where $c \in \{0, 1, \ldots, p-1\}$ is a residue modulo $p$ and $\varepsilon = \pm$ is a sign. Thus $2+$ and $3-$ are examples of 5-signatures.

A $(p)$-*signature sequence* $X$ is a sequence

$$X : \quad c_1\varepsilon_1 \, c_2\varepsilon_2 \, \cdots \, c_s\varepsilon_s$$

where each $c_i\varepsilon_i$ is a signature.

For $0 \le i \le s$ and $0 \le \alpha \le p-1$ we define

$$\sigma_\alpha^X(i) = \sum_{\{k \le i | c_k = \alpha\}} \varepsilon_k \,.$$

We make the conventions that an empty sum is 0 and that $+$ is counted as $+1$ and $-$ as $-1$ in the sum. The *end value* $\sigma_\alpha^X$ of $\alpha$ in $X$ is then defined to be

$$\sigma_\alpha^X = \sigma_\alpha^X(s) \,.$$

Let

$$\pi_\alpha(X) = \max\{\sigma_\alpha^X(i) \mid 0 \le i \le s\} \,.$$

We call $c_i$ *normal* of residue $\alpha$ if $\sigma_\alpha^X(i) > \sigma_\alpha^X(j)$ for all $j \le i-1$ and $\sigma_\alpha^X(i) > 0$. This is only possible when $c_i\varepsilon_i = \alpha+$. In this case we also call $\sigma_\alpha^X(i)$ the *height* ht $c_i$ of $c_i$. Moreover, $c_i$ is called *good* of residue $\alpha$ if $c_i$ is normal of residue $\alpha$ and $i$ is minimal with

$$\sigma_\alpha^X(i) = \pi_\alpha(X) \,.$$

Note that if $c_i$ is good of residue $\alpha$ then ht $c_i = \pi_\alpha(X)$.

The *node sequence* $N(\lambda)$ of $\lambda$ is the signature sequence consisting of the residues of the indent and removable nodes of $\lambda$, read from left to right, top to bottom in $\lambda$. For each indent residue the sign is $+$ and for each removable residue the sign is $-$.

**Example.** Let $p = 5$ and $\lambda = (12, 7^2, 5^3, 3, 1^3)$ as above. Then the node sequence of $\lambda$ is

$$N(\lambda): \quad \underline{1+} \;\; 2- \;\; 1- \;\; \underline{4+} \;\; 2- \;\; \underline{4+} \;\; \underline{1+} \;\; 2- \;\; 4- \;\; \underline{1+} \;\; 0-$$

where we have underlined the normal entries. *Note that these entries and their heights correspond to normal nodes of the same height in $Y(\lambda)$.*

These combinatorial concepts play a vital rôle in the following theorem due to Kleshchev (see [6]-[9]).

**Theorem 2.1** *Let $\lambda$ be a $p$-regular partition of $n$, $n \in \mathbb{N}$, $n \ge 2$. Then the following holds:*

*(i)* $soc\,(D^\lambda|_{S_{n-1}}) \simeq \bigoplus_{A \text{ good}} D^{\lambda \backslash A}$.

*(ii)* $D^\lambda|_{S_{n-1}}$ *is completely reducible if and only if all normal nodes in $\lambda$ are good.*

*(iii) Let $A$ be a removable node of $\lambda$ such that $\lambda \backslash A$ is $p$-regular. Then the multiplicity of $D^{\lambda\backslash A}$ in $D^\lambda|_{S_{n-1}}$ is given by*

$$\left[D^\lambda|_{S_{n-1}} : D^{\lambda\backslash A}\right] = \begin{cases} \text{ht } A & \text{if } A \text{ is normal in } \lambda \\ 0 & \text{else} \end{cases}$$

In particular, this theorem tells us for which $p$-regular partitions $\lambda$ of $n$ we have that the restriction $D^\lambda|_{S_{n-1}}$ is irreducible. Indeed, by (ii) all normal nodes in $\lambda$ have to be good and by (i) there can only be one good node in $\lambda$. Now obviously the removable node in the first block of $\lambda$ is normal. Since the first terms of $N(\lambda)$ are

$$l_1 - \ (l_1 - a_1) + \ (l_2 - a_1) - \ (l_2 - a_1 - a_2) +$$

we see that in order that the second removable node is not normal we need that $l_1 \equiv l_2 - a_1 - a_2$, i.e. $\beta(1,2) \equiv 0$. Continuing in this way we obtain that the conditions

$$\beta(i, i+1) \equiv 0 (\operatorname{mod} p)$$

for $1 \leq i \leq t$ are necessary and sufficient for $\lambda$ to have exactly one normal node. This leads us to the following definition.

**Definition 2.2** *Let $\lambda = (l_1^{a_1}, \ldots, l_t^{a_t})$ be a $p$-regular partition, where $l_1 > l_2 > \cdots > l_t$, $0 < a_i < p$ for $i = 1, \ldots, t$. Then $\lambda$ is called a* JS-partition *if its parts satisfy the congruence*

$$\beta(i, i+1) = l_i - l_{i+1} + a_i + a_{i+1} \equiv 0 \bmod p \quad \textit{for } 1 \leq i < t\ .$$

*The* type *$\alpha$ of $\lambda$ is defined to be the residue $\alpha$ of $l_1 - a_1$ mod $p$.*

Note that the type of a JS-partition is just the residue of the unique normal (and thus good) node. The letters JS are an abbreviation of Jantzen-Seitz. These authors conjectured the equivalence of $(i) \Leftrightarrow (iii)$ in the following consequence of Kleshchev's theorem. With the above, we have proved the equivalence $(ii) \Leftrightarrow (iii)$ below.

**Corollary 2.3** *With notation as above, the following are equivalent:*

*(i) $D^\lambda|_{S_{n-1}}$ is irreducible.*

*(ii) $\lambda$ has exactly one normal node (which is then the only good node in $\lambda$).*

*(iii) $\lambda$ is a JS-partition.*

Next we want to introduce a different description of $p$-regular partitions which originated from the definition of the Mullineux map as given by Mullineux [15].

Let $\lambda$ be a $p$-regular partition of $n$. The *$p$-rim* of $\lambda$ is a part of the rim of $\lambda$ ([7], p. 56), which is composed of *$p$-segments*. Each $p$-segment except possibly the last contains $p$ points. The first $p$-segment consists of the first $p$ points of the rim of $\lambda$, starting with the longest row. (If the rim contains at most $p$ points it is the entire rim.) The next segment is obtained by starting in the row next below the previous $p$-segment. This process is continued until the final row is reached. We let $a_1$ be the number of nodes in the $p$-rim of $\lambda = \lambda^{(1)}$ and let $r_1$ be the number of rows in $\lambda$. Removing the $p$-rim of $\lambda = \lambda^{(1)}$ we get a new $p$-regular partition $\lambda^{(2)}$ of $n - a_1$.

We let $a_2, r_2$ be the length of the $p$-rim and the number of parts of $\lambda^{(2)}$ respectively. Continuing in this way we get a sequence of partitions $\lambda = \lambda^{(1)}, \lambda^{(2)}, \ldots, \lambda^{(m)}$ where $\lambda^{(m)} \neq 0$ and $\lambda^{(m+1)} = 0$; the corresponding *Mullineux symbol* of $\lambda$ is then given by

$$G_p(\lambda) = \begin{pmatrix} a_1 & a_2 & \cdots & a_m \\ r_1 & r_2 & \cdots & r_m \end{pmatrix} .$$

The integer $m$ is called the *length* of the symbol.

If $p|a_i$, we call the corresponding column $\begin{smallmatrix} a_i \\ r_i \end{smallmatrix}$ of the Mullineux symbol a *singular* column, otherwise the column is called *regular*.

For some properties of the Mullineux symbol we refer to [1], [2], [15], [16].

It turns out that it is often preferable to work with the *Residue symbol* of $\lambda$ introduced in [2] rather than with the Mullineux symbol.

If the Mullineux symbol $G_p(\lambda)$ is given as above, then the *Residue symbol* $R_p(\lambda)$ of $\lambda$ is defined as

$$R_p(\lambda) = \left\{ \begin{matrix} x_1 & x_2 & \cdots & x_m \\ y_1 & y_2 & \cdots & y_m \end{matrix} \right\}$$

where $x_j$ is the residue of $a_{m+1-j} - r_{m+1-j}$ modulo $p$ and $y_j$ is the residue of $1 - r_{m+1-j}$ modulo $p$. The Mullineux symbol $G_p(\lambda)$ can be recovered from the Residue symbol $R_p(\lambda)$ because the entries in the Mullineux symbol satisfy strong restrictions [1].

Now we turn to the *Mullineux* (signature) *sequence* $M(\lambda)$, which is defined as follows:

Let the residue symbol of $\lambda$ be as above. Then

$$\begin{array}{rlllll} M(\lambda) = 0- & x_1+ & (x_1+1)- & y_1+ & (y_1-1)- \\ & x_2+ & (x_2+1)- & y_2+ & (y_2-1)- \\ & & \vdots & \vdots & \\ & x_m+ & (x_m+1)- & y_m+ & (y_m-1)- \end{array}$$

Starting with the signature $0-$ corresponds to starting with an empty partition at the beginning which just has the indent node $(1,1)$ of residue 0.

In [2] the following result was proved.

**Theorem 2.4** *Let $\lambda$ be a $p$-regular partition.*
*Then for all $\alpha$, $0 \leq \alpha \leq p-1$ we have*

$$\pi_\alpha(M(\lambda)) = \pi_\alpha(N(\lambda)) .$$

From the remarks above we then obtain

**Corollary 2.5** *The following statements are equivalent for a $p$-regular partition $\lambda$ $(0 \leq \alpha \leq p-1)$.*

*(i)* $\lambda$ *has a normal (good) node of residue* $\alpha$*.*

*(ii)* $M(\lambda)$ *has a normal (good) entry of residue* $\alpha$*.*

*(iii)* $N(\lambda)$ *has a normal (good) entry of residue* $\alpha$*.*

# 3 Residue symbols and cores for JS-partitions

Based on the results 2.3 and 2.5 in the previous section we present here a simplified proof of the following characterisation of the residue symbols for JS-partitions. The construction rules in the theorem below may also be obtained by translating the main result of [1] to residue symbols but with the proof below we can considerably shorten the lengthy and tedious case-by-case analysis of the previous proof. This illustrates once again the usefulness of the residue symbols.

For a $p$-regular partition $\lambda$ with Mullineux sequence $X = M(\lambda)$, we put $\sigma_\alpha = \sigma_\alpha^X$ for all residues $\alpha \in \{0, \ldots, p-1\}$ and define the *end value vector* $\sigma$ of $\lambda$ by $\sigma = (\sigma_0, \sigma_1, \ldots, \sigma_{p-1})$.

Below, $e_i$ denotes the vector of length $p$ with 1 at residue position $i$, 0 elsewhere.

**Theorem 3.1** *The residue symbols of JS-partitions of type* $\alpha$ *can be constructed iteratively as described below.*

*The residue symbols* $R_p(\lambda) = \left\{ \begin{matrix} x_1 \\ y_1 \end{matrix} \right\}$ *of length* 1 *of JS-partitions* $\lambda$ *of type* $\alpha$ *are*

$$\left\{ \begin{matrix} x_1 \\ y_1 \end{matrix} \right\} = \left\{ \begin{matrix} 0 \\ \alpha \end{matrix} \right\} (\textit{with } \alpha \neq 1) \quad \textit{or} \quad \left\{ \begin{matrix} \alpha \\ 0 \end{matrix} \right\} \quad \textit{or} \quad \left\{ \begin{matrix} \alpha \\ \alpha+1 \end{matrix} \right\} (\textit{with } \alpha \neq 0).$$

*If* $\lambda$ *is a JS-partition of type* $\alpha$ *with residue symbol* $R_p(\lambda) = \left\{ \begin{matrix} x_1 & \cdots & x_k \\ y_1 & \cdots & y_k \end{matrix} \right\}$*, then we have the following possibilities for extending the residue symbol of* $\lambda$ *to a residue symbol of a JS-partition* $R_p(\mu) = \left\{ \begin{matrix} x_1 & \cdots & x_k & x_{k+1} \\ y_1 & \cdots & y_k & y_{k+1} \end{matrix} \right\}$ *of type* $\alpha$ *of length* $k+1$*:*

$$\begin{matrix} x_{k+1} \\ y_{k+1} \end{matrix} = \begin{cases} \begin{matrix} \alpha+1-y_k \\ y_k - 1 \end{matrix} & (\textit{if } 2y_k \not\equiv \alpha+3) \\ \quad \textit{or} & \qquad\qquad (\textit{regular extensions}) \\ \begin{matrix} y_k - 1 \\ \alpha+1-y_k \end{matrix} & (\textit{if } 2y_k \not\equiv \alpha+1) \\ \quad \textit{or} & \\ \begin{matrix} y_k - 1 \\ y_k \end{matrix} & \\ \quad \textit{or} & \qquad\qquad (\textit{singular extensions}) \\ \begin{matrix} \alpha+1-y_k \\ \alpha+2-y_k \end{matrix} & \end{cases}$$

*Moreover, one of the following holds for the end value vector* $\sigma$ *of* $R_p(\lambda)$ *and the entry* $y_k$*:*

*(a)* $\sigma = v_0 = -e_0$ *and* $y_k \in \{1, \alpha+1\}$.

*(b)* $\sigma = v_{\alpha,\beta} = e_\alpha - e_\beta - e_{\alpha-\beta}$, $\alpha \neq \beta$, *and* $y_k \in \{1+\beta, \alpha+1-\beta\}$.

*(c)*

$$\sigma = w_\alpha = \begin{cases} e_\alpha - 2e_{\alpha/2} & \text{if } \alpha \neq 0 \text{ even, and in this case } y_k = \frac{\alpha}{2}+1 \\ e_\alpha - 2e_{(p+\alpha)/2} & \text{if } \alpha + p \text{ even, and in this case } y_k = \frac{\alpha+p}{2}+1 \end{cases}.$$

**Proof.** By 2.5 and 2.3, $\lambda$ is a JS-partition of type $\alpha$ if and only if its Mullineux sequence $X$ satisfies

$$\sigma^X_\beta(i) \leq \begin{cases} 0 & \text{if } \beta \neq \alpha \\ 1 & \text{if } \beta = \alpha \end{cases}.$$

for all $i$ and all residues $\beta$. This leads to very strong conditions on the columns in the residue symbol.

First we consider the case of residue symbols $R_p(\lambda) = \left\{ \begin{matrix} x_1 \\ y_1 \end{matrix} \right\}$ of JS-partitions $\lambda$ of type $\alpha$ which are of length $k = 1$.

Here, the Mullineux sequence is

$$M(\lambda) = 0- \quad x_1+ \quad (x_1+1)- \quad y_1+ \quad (y_1-1)-$$

By the above characterization, it is then clear that the only possible residue symbols are

$$\left\{ \begin{matrix} 0 \\ \alpha \end{matrix} \right\} \text{(with } \alpha \neq 1\text{)}, \left\{ \begin{matrix} \alpha \\ 0 \end{matrix} \right\} \text{ and } \left\{ \begin{matrix} \alpha \\ \alpha+1 \end{matrix} \right\} \text{(with } \alpha \neq 0\text{)}$$

as claimed (note that $\left\{ \begin{matrix} 0 \\ 1 \end{matrix} \right\}$ is not the residue symbol of a $p$-regular partition).

For these three symbols we have the following end value vectors. For $\left\{ \begin{matrix} 0 \\ \alpha \end{matrix} \right\}$, we have

$$\sigma = v_{\alpha,1} \quad (= w_2 \text{ if } \alpha = 2),$$

and $y_1 = \alpha$ satisfies the conditions stated in the Theorem.

For $\left\{ \begin{matrix} \alpha \\ 0 \end{matrix} \right\}$,

$$\sigma = \begin{cases} v_{\alpha,p-1} & \\ w_{p-2} & \text{if } \alpha = p-2 \end{cases}$$

and $y_1 = 0$ satisfies the required conditions.

Finally, for $\left\{ \begin{matrix} \alpha \\ \alpha+1 \end{matrix} \right\}$, $\sigma = v_0$, and $y_1 = \alpha + 1$ is as required.

Now assume that we have already proved for a given $k \in \mathbb{N}$ that any residue symbol of a JS-partition $\lambda$ of type $\alpha$ of length $k$ as above satisfies the stated assertions on its end value vector and its entry $y_k$. We now investigate the possible extensions of $R_p(\lambda)$ to a residue symbol $R_p(\mu) = \left\{ \begin{matrix} x_1 & \cdots & x_k & x_{k+1} \\ y_1 & \cdots & y_k & y_{k+1} \end{matrix} \right\}$ of type $\alpha$ of length $k+1$.

In case $(a)$, when $\sigma = v_0$ and $y_k \in \{1, \alpha+1\}$, we have the following possibilities:

$$\begin{matrix} x_{k+1} \\ y_{k+1} \end{matrix} = \begin{matrix} 0 \\ \alpha \end{matrix} , \begin{matrix} \alpha \\ 0 \end{matrix} , \begin{matrix} \alpha \\ \alpha+1 \end{matrix} \text{ or } \begin{matrix} 0 \\ 1 \end{matrix} .$$

These are exactly the extensions we obtain by the construction rules given in the Theorem, using that $y_k = 1$ or $y_k = \alpha + 1$.

Furthermore, in these cases, the new end value vectors are

$$v_{\alpha,1} , v_{\alpha,p-1} , v_0 , v_0$$

respectively, and since $y_{k+1} = \alpha , 0 , \alpha+1 , 1$ in the four cases, the claimed restrictions are satisfied for the extended residue symbol.

In case $(b)$, when $\sigma = v_{\alpha,\beta} = e_\alpha - e_\beta - e_{\alpha-\beta}$, $\alpha \neq \beta$, and $y_k \in \{1+\beta, \alpha+1-\beta\}$, we have the following possibilities:

$$\begin{matrix} x_{k+1} \\ y_{k+1} \end{matrix} = \begin{matrix} \beta \\ \alpha-\beta \end{matrix} , \begin{matrix} \alpha-\beta \\ \beta \end{matrix} , \begin{matrix} \beta \\ \beta+1 \end{matrix} \text{ or } \begin{matrix} \alpha-\beta \\ \alpha-\beta+1 \end{matrix} .$$

Again, the construction rules in the Theorem, applied to $y_k = 1+\beta$ or $y_k = \alpha+1-\beta$ give the same symbols.

In these four cases, the new end value vectors are

$$v_{\alpha,\beta-1} , v_{\alpha,\beta+1} , v_{\alpha,\beta} , v_{\alpha,\beta}$$

and one immediately checks that the corresponding $y_{k+1}$ satisfy the given conditions.

In case $(c)$, we consider the situation when $\sigma = w_\alpha$, $\alpha \neq 0$ even, and $y_k = \frac{\alpha}{2} + 1$; the other case in $(c)$ is completely analogous. Here we only have the following two possibilities:

$$\begin{matrix} x_{k+1} \\ y_{k+1} \end{matrix} = \begin{matrix} \alpha/2 \\ \alpha/2 \end{matrix} \text{ or } \begin{matrix} \alpha/2 \\ \alpha/2+1 \end{matrix}$$

Again, the construction rules in the Theorem, applied to $y_k = \frac{\alpha}{2} + 1$ give the same symbols. The extended residue symbol then has the end value vector $v_{\alpha,\alpha/2+1}$ resp. $w_\alpha$, and $y_{k+1}$ satisfies the corresponding restrictions.

Since the argument given in the first paragraph of the proof implies that the beginning of a residue symbol of a JS-partition of type $\alpha$ is itself a residue symbol of a JS-partition of type $\alpha$, any such residue symbol is constructed iteratively from a residue symbol of a JS-partition of type $\alpha$ of length 1.

Thus all the assertions of the Theorem are proved. $\diamond$

We note the following consequence of the construction rules for JS-partitions which allows to get rid of the singular columns in the residue symbol for some computations.

**Corollary 3.2** *Let* $R_p(\lambda) = \left\{ \begin{matrix} x_1 & x_2 & \cdots & x_m \\ y_1 & y_2 & \cdots & y_m \end{matrix} \right\}$ *be a JS-partition of type* $\alpha$*. If there are only singular columns between the regular columns* $\substack{x_l \\ y_l}$ *and* $\substack{x_k \\ y_k}$ *(where* $l < k$*), then* $\substack{x_k \\ y_k} = \substack{\alpha+1-y_l \\ y_l-1} = \substack{x_l+1 \\ \alpha-(x_l+1)}$ *or* $\substack{y_l-1 \\ \alpha+1-y_l}$.

*Also, if* $R_p(\lambda)$ *starts with singular columns, the first regular column is of the form* $\substack{\alpha \\ 0}$ *or* $\substack{0 \\ \alpha}$.

**Proof.** For $l = k-1$ the assertion certainly holds.

Assume now that $l < k-1$.

The singular extensions of $\substack{x_l \\ y_l}$ are $\substack{y_l-1 \\ y_l}$ and $\substack{\alpha+1-y_l \\ \alpha+2-y_l}$. Extending these singularly gives again the same two possible columns. So this repeats until column $k-1$.

From these two columns we then get in both cases the regular extensions $\substack{y_l-1 \\ \alpha+1-y_l}$ and $\substack{\alpha+1-y_l \\ y_l-1}$.

The second statement is proved similarly. $\diamond$

The construction rules given in Theorem 3.1 can very nicely be described by a suitable diagram; this also explains the Corollary above in a very natural way. We formulate this alternative description (which is an immediate translation of the rules given in the previous Theorem) as

**Theorem 3.3** *The residue symbols of JS-partitions of type* $\alpha$ *can be iteratively constructed along the directed graph given below in the following way.*

*Start at the vertices of the graph labeled by*

$$\left\{ \begin{matrix} 0 \\ \alpha \end{matrix} \right\} (\textit{if } \alpha \neq 1) \quad \textit{or} \quad \left\{ \begin{matrix} \alpha \\ 0 \end{matrix} \right\} \quad \textit{or} \quad \left\{ \begin{matrix} \alpha \\ \alpha+1 \end{matrix} \right\} (\textit{if } \alpha \neq 0);$$

*this is the first column in the residue symbol to be constructed. Then add on columns for the residue symbol by following a path in the graph. In the diagrams below we omit the brackets for simplification; the symbol* $\cap$ *means that we have a loop at the corresponding vertex of the graph. For simplification, we have drawn a double-vertex in the middle of each square (and on one side) which is labeled by two singular columns; the residue symbol can pick up any of the two columns, and in looping around the vertex again any of the two columns can be chosen (so the two arrows going into (resp. out) of the double-vertex correspond to four arrows altogether, and the loop in the middle corresponds to two loops and a pair of opposite arrows between the two singular columns).*

*For $\alpha = 0$ the graph is*

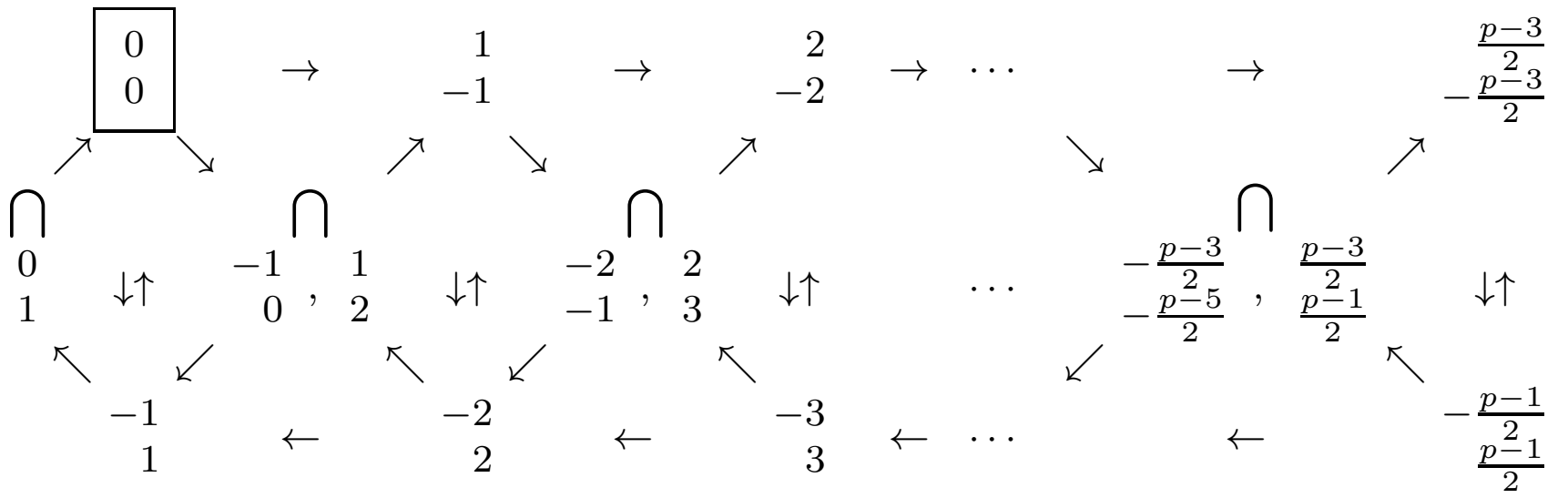

*For $\alpha \neq 0$ even the graph is as shown in Figure 1 (the graph for $\alpha$ odd is analogous). In these graphs, we have marked the admissible initial columns.*

We next want to compute the $p$-cores of JS-partitions. We recall that the *$p$-core* $\lambda_{(p)}$ of a partition $\lambda$ is obtained by removing $p$-hooks as long as possible; while the removal process is not unique the resulting $p$-regular partition is unique (we refer to [7] or [18] for a more detailed introduction into this notion and its properties). We define the *weight* $w$ of $\lambda$ by $w = (|\lambda| - |\lambda_{(p)}|)/p$.

The $p$-core of $\lambda$ determines the $p$-block to which an ordinary or modular irreducible character labeled by $\lambda$ belongs. The *weight* of a $p$-block is the common weight of the partitions labeling the characters in the block.

The *$p$-content* $c(\lambda) = (c_0, \ldots, c_{p-1})$ of a partition $\lambda$ is defined by counting the number of nodes of a given residue in the $p$-residue diagram of $\lambda$, i.e. $c_i$ is the number of nodes of $\lambda$ of $p$-residue $i$. In the example given before, the 5-content of $\lambda = (6^2, 5, 4)$ is $c(\lambda) = (c_0, \ldots, c_4) = (5, 3, 4, 4, 5)$.

We recall that the $p$-content determines the $p$-core of a partition. To achieve this explicitly, for given $c = (c_0, c_1, \ldots, c_{p-1})$ we define the associated $\vec{n}$-vector by $\vec{n} = (c_0 - c_1, c_1 - c_2, \ldots, c_{p-2} - c_{p-1}, c_{p-1} - c_0)$. Now, for any vector

$$\vec{n} \in \{(n_0, \ldots, n_{p-1}) \in \mathbb{Z}^p \mid \sum_{i=0}^{p-1} n_i = 0\}$$

there is a unique $p$-core $\mu$ with this $\vec{n}$-vector $\vec{n}$ associated to its $p$-content $c(\mu)$ (for short, we also say that $\vec{n}$ is the $\vec{n}$-vector associated to $\mu$). We refer the reader to [6] for the description of this bijection.

In [2] we have shown how to find the $\vec{n}$-vector via the Mullineux or residue symbol; this is based on the fact that the contribution of a regular column $\begin{smallmatrix}\alpha\\ \beta\end{smallmatrix}$ in the residue symbol of $\lambda$ for the $\vec{n}$-vector of $\lambda$ corresponds to that of a hook with arm node of residue $\alpha$ and leg node of residue $\beta$, whereas a singular column gives no contribution.

**Proposition 3.4** *[2] Let $\lambda$ be a $p$-regular partition with residue symbol*

$$R_p(\lambda) = \left\{ \begin{matrix} x_1 & \cdots & x_m \\ y_1 & \cdots & y_m \end{matrix} \right\}.$$

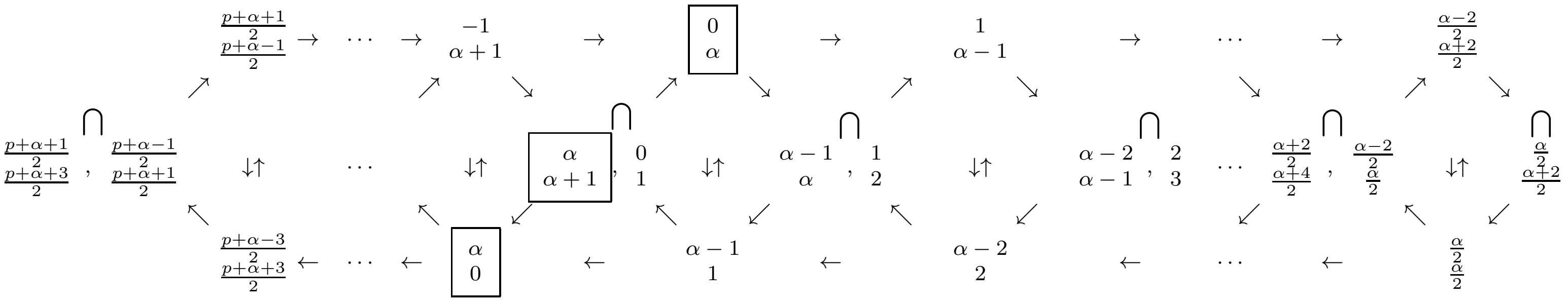

**Figure 1 : The construction diagram for JS-partitions of type $\alpha$**

*Then the associated $\vec{n}$-vector $\vec{n} = (n_0, \ldots, n_{p-1})$ is given by*

$$n_j = |\{i \mid x_i = j\}| - |\{i \mid y_i = j+1\}|$$

By comparing the $\vec{n}$-vectors we will now prove that JS-partitions of type $\alpha$ have a rectangular shaped $p$-core of type $\alpha$ or an empty $p$-core, and we determine this $p$-core precisely from the final column of the residue symbol.

**Theorem 3.5** *If $\lambda$ is a JS-partition of type $\alpha$, then its $p$-core is either empty or is rectangular shaped with removable node of residue $\alpha$, i.e. $\lambda_{(p)} = \emptyset$ or $\lambda_{(p)} = (l^a)$ for some $l, a \in \mathbb{N}$ with $l - a \equiv \alpha$.*

*More precisely, if*

$$R_p(\lambda) = \left\{ \begin{matrix} x_1 & \cdots & x_m \\ y_1 & \cdots & y_m \end{matrix} \right\},$$

*then*

$$\lambda_{(p)} = \begin{cases} (l^{p+l-\alpha}) & \text{if } \begin{matrix} x_m \\ y_m \end{matrix} \in \left\{ \begin{matrix} l-1 \\ \alpha-l+1 \end{matrix}, \begin{matrix} \alpha-l-1 \\ l+1 \end{matrix}, \begin{matrix} l \\ l+1 \end{matrix}, \begin{matrix} \alpha-l \\ \alpha-l+1 \end{matrix} \right\} \\ & \quad \text{and } 1 \leq l \leq \frac{\alpha}{2} \\ (l^{l-\alpha}) & \text{if } \begin{matrix} x_m \\ y_m \end{matrix} \in \left\{ \begin{matrix} l-1 \\ \alpha-l+1 \end{matrix}, \begin{matrix} \alpha-l-1 \\ l+1 \end{matrix}, \begin{matrix} l \\ l+1 \end{matrix}, \begin{matrix} \alpha-l \\ \alpha-l+1 \end{matrix} \right\} \\ & \quad \text{and } \alpha+1 \leq l < \frac{p+1+\alpha}{2} \\ \emptyset & \text{if } \begin{matrix} x_m \\ y_m \end{matrix} \in \left\{ \begin{matrix} \alpha-1 \\ 1 \end{matrix}, \begin{matrix} -1 \\ \alpha+1 \end{matrix}, \begin{matrix} 0 \\ 1 \end{matrix}, \begin{matrix} \alpha \\ \alpha+1 \end{matrix} \right\} \end{cases}$$

**Remark.** Note that the 'missing' values for $l$ cannot occur because they would not give a $p$-core.

**Proof.** Let $\mu = (l^a)$ be a rectangular partition; $\mu$ is a $p$-core if and only if $l + a - 1 < p$. The removable corner of $\mu$ is of residue $\alpha = l - a$.

We compute the $\vec{n}$-vector of $\mu$.

*Case 1.* First assume $a \leq l$. Then the residue diagram of $\mu$ is

$$\boxed{\begin{array}{cccccccc} 0 & 1 & 2 & \cdots & \alpha & \cdots & & l-1 \\ p-1 & 0 & 1 & 2 & \cdots & \alpha & & \vdots \\ \vdots & & \ddots & & & & \ddots & \\ p-a+1 & \cdots & & 0 & \cdots & & & \alpha \end{array}}$$

So the $\vec{n}$-vector of $\mu$ is

$$\begin{array}{rcllllllll} \vec{n}(\mu) & = & (0,\dots,0, & 1, & \dots, & 1, & 0,\dots,0, & -1, & \dots, & -1\,) \\ & & & l-a & & l-1 & & p-a & & p-1 \\ & & & =\alpha & & & & & & \\ & =: & b_a^{\alpha+} & & & & & & & \end{array}$$

Here we have marked the residues $i$ for the relevant entries $n_i$ of the $\vec{n}$-vector.

Set $b_0 = b_0^{\alpha+} = b_0^{\alpha-} = 0$.

*Case 2.* Now assume $l < a$. For any partition $\nu$ with $\vec{n}$-vector $\vec{n}(\nu) = (n_0, \dots, n_{p-1})$, the $\vec{n}$-vector of the conjugate partition $\nu'$ is given by $\vec{n}(\nu') = (-n_{p-1}, \dots, n_0)$.

Thus we can apply the previous case and the $\vec{n}$-vector of $\mu$ is

$$\begin{array}{rcllllllll} \vec{n}(\mu) & = & (1, & \dots, & 1, & 0,\dots,0, & -1, & \dots, & -1, & 0,\dots,0) \\ & & 0 & & l-1 & & p-a & & \alpha-1 & \\ & & & & & & =\alpha-l & & & \\ & =: & b_l^{\alpha-} & & & & & & & \end{array}$$

Now it suffices to prove that a JS-partition $\lambda$ of type $\alpha$ has an $\vec{n}$-vector of the form above, with the predicted dependency on its final column.

From a column $\left\{\begin{smallmatrix} x \\ y \end{smallmatrix}\right\}$ in $R_p(\lambda)$ we get the contribution 1 to $n_x$ and $-1$ to $n_{y-1}$ and nothing else (note that for $y = x+1$ the total contribution to $n_x = n_{y-1}$ is also 0).

For $m = 1$, we have

$$\vec{n}(\left\{\begin{matrix} 0 \\ \alpha \end{matrix}\right\}) = b_1^{\alpha-} = \vec{n}((1^a)) \quad \text{with } a = p+1-\alpha$$

and

$$\vec{n}(\left\{\begin{matrix} \alpha \\ 0 \end{matrix}\right\}) = b_1^{\alpha+} = \vec{n}((l)) \quad \text{with } l = \alpha+1$$

and

$$\vec{n}(\left\{\begin{matrix} \alpha \\ \alpha+1 \end{matrix}\right\}) = 0 = \vec{n}(\emptyset)$$

If the residue symbol has only singular columns, then we had started at $\begin{smallmatrix} \alpha \\ \alpha+1 \end{smallmatrix}$ and have only added columns of the same type or of the form $\begin{smallmatrix} 0 \\ 1 \end{smallmatrix}$, and the $p$-core is clearly empty.

So we may assume that the residue symbol has a regular column; as the regular predecessors of the singular columns given in the sets for the three cases above all are regular columns in the same set, we may assume that the final column in the residue symbol is regular.

By Corollary 3.2 and the fact that singular columns do not contribute to the $\vec{n}$-vector, we may then assume that our residue symbol is constructed solely by regular extensions (and from a regular start). In particular, we have $x_i + y_i = \alpha \pmod p$ for all $i \in \{1, \dots, m\}$.

Let $\vec{n}_i$ be the $\vec{n}$-vector of the residue symbol up to column $i$. For $j \in \{0, \dots, p-1\}$ set $v_j^\alpha = e_j - e_{\alpha-1-j}$, where $e_j$ denotes the vector with 1 at residue position $j$, 0 elsewhere. (Perhaps the best way of understanding the change of the $\vec{n}$-vector and thus the $p$-core along the construction of the residue symbol in the graph given before is to write these vectors $v_j^\alpha = e_j - e_{\alpha-1-j}$ at the corresponding vertices $\substack{j \\ \alpha-j}$ in the graph and watch what happens along the path; remember that there are no contributions to the $\vec{n}$-vector from the singular columns.)

Now assume that the statement has already been proved up to column $m-1$. The regular predecessors of $\substack{x_m \\ y_m}$ are of the form $\substack{x_m-1 \\ y_m+1}$ and $\substack{y_m-1 \\ x_m+1}$, and in both cases we have the same $p$-core resp. $\vec{n}$-vector $\vec{n}_{m-1}$ attached to the penultimate step for the residue symbol.

In the first of the three cases in the statement of the Theorem this is (by induction and the description given at the beginning of this proof)

$$\vec{n}_{m-1} = b_{x_m}^{\alpha-} = b_{l-1}^{\alpha-}$$

and the contribution of the column $\substack{x_m \\ y_m}$ then gives

$$\vec{n}_m = b_l^{\alpha-}$$

as desired.

The other cases are completely analogous. $\diamond$

Alternatively, the theorem above can also easily be proved by computing the $\vec{n}$-vector using the orginal definition of JS-partitions. We carry this out and state explicitly how the $p$-core depends on the length of a JS-partition of type $\alpha$:

**Theorem 3.6** *If $\lambda$ is a JS-partition of type $\alpha$ and of length $l \equiv r \in \{0, \dots, p-1\}$ (mod p), then its p-core is*

$$\lambda_{(p)} = \begin{cases} ((\alpha+s)^s) & \text{with} \quad s = r & \text{if } r \le \frac{p-\alpha}{2} \\ & \qquad\quad s = (p-\alpha) - r & \text{if } \frac{p-\alpha}{2} < r \le p-\alpha \\ (s^{p-\alpha+s}) & \text{with} \quad s = r - (p-\alpha) & \text{if } p-\alpha < r \le p - \frac{\alpha}{2} \\ & \qquad\quad s = p - r & \text{if } p - \frac{\alpha}{2} < r < p \end{cases}$$

**Proof.** We compute the $\vec{n}$-vector of $\lambda$ by considering the contributions of each row. Counting down the beginnings of the rows, we get $l$ contributions $-1$, for $p-1$ down to $p-l$. At the end of the rows, we start with the good corner node $A$ of residue $\alpha$ and go up to the end node $B$ of the first row, which is of residue $\beta$ (say). By the JS-condition, the second corner node is then of residue $\beta+1$; hence in jumping from $B$ to this node and going up the second block, we are just obtaining contributions 1 for the successive sequence of residues $\alpha$, $\alpha+1$, $\alpha+2$, .... Continuing in this way, we clearly have exactly $l$ contributions 1 for the successive residues $\alpha$, ..., $\alpha+l-1$.

Most of these cancel out against the $-1$ contributions from the beginnings of the rows and we are left with the following possible $\vec{n}$-vectors:

$$\vec{n}(\lambda) = \begin{cases} b_s^{\alpha+} & \text{with} \quad s = r & \text{if } r \leq \frac{p-\alpha}{2} \\ & \qquad\quad s = (p-\alpha) - r & \text{if } \frac{p-\alpha}{2} < r \leq p - \alpha \\ b_s^{\alpha-} & \text{with} \quad s = r - (p-\alpha) & \text{if } p - \alpha < r \leq p - \frac{\alpha}{2} \\ & \qquad\quad s = p - r & \text{if } p - \frac{\alpha}{2} < r < p \end{cases}$$

In the previous proof we had already computed the $\vec{n}$-vectors for rectangles of type $\alpha$, and thus we deduce the statement in the Theorem immediately from the description above. $\diamond$

**Remark.** The fact that the $p$-cores of JS-partitions are empty or of rectangular shape has also been proved recently by algebraic methods by Foda et al. [4, Theorem 6.3] in the context of studying certain exactly solvable models in statistical mechanics and their relations to Hecke algebras.

Finally, we want to consider the question of existence of JS-partitions of a given rectangular (or empty) $p$-core and given weight $w$.

**Proposition 3.7** *For any $p$-core of the form $\mu = \emptyset$ or $\mu = (l^a)$ and for any given weight $w$ there exists a JS-partition with this $p$-core $\mu$ and weight $w$.*

**Proof.** For $\mu = \emptyset$, obviously $\lambda = (pw)$ is a JS-partition of weight $w$ with $\lambda_{(p)} = \emptyset$. For a $p$-core $\mu = (l^a)$, its Mullineux symbol is

$$G_p(\mu) = \begin{pmatrix} l+a-1 & l+a-3 & \cdots & l-a+1 \\ a & a-1 & \cdots & 1 \end{pmatrix}$$

where $l + a - 1 < p$ since $\mu$ is a $p$-core. Hence we can extend this Mullineux symbol by $w$ columns $\begin{smallmatrix} p \\ a \end{smallmatrix}$ and obtain the Mullineux symbol of a JS-partition with $p$-core $\mu$ and of weight $w$. $\diamond$

For $a \in \mathbb{N}$, we define the *$p$-level* of $a$ to be $\left[\frac{a}{p}\right]$.

**Proposition 3.8** *Let* $\begin{smallmatrix} x' \\ y' \end{smallmatrix} \rightarrow \begin{smallmatrix} x \\ y \end{smallmatrix}$ *be an edge in the construction diagram of JS-partitions of type $\alpha$. Define $t, t', s \in \{0, 1, \ldots, p-1\}$ by*

$$t \equiv x - y + 1\,,\ t' \equiv x' - y' + 1\,,\ s \equiv y' - y \quad (mod\ p)\,,$$

*and let*

$$\varepsilon = \begin{cases} 0 & if\ t = 0 \\ 1 & otherwise \end{cases}\,,\ \varepsilon' = \begin{cases} 0 & if\ t' = 0 \\ 1 & otherwise \end{cases}.$$

*Then the increase $d$ of the $p$-level of the length of the $p$-rim in extending the residue symbol along the given edge is independent of the previous construction steps of the residue symbol, and it is given as follows:*

$$d = \begin{cases} [\frac{t'+2s}{p}] & \text{if } \varepsilon' = \varepsilon = 1 \\ [\frac{t'+s+p}{p}] & \text{if } \varepsilon' = 1,\, \varepsilon = 0 \\ [\frac{p-1+s-t}{p}] & \text{if } \varepsilon' = 0 \end{cases} .$$

**Proof.** Let $\substack{a' \\ r'}$ and $\substack{a \\ r}$ be the columns in the Mullineux symbols corresponding to $\substack{x' \\ y'}$ and $\substack{x \\ y}$, respectively. By the definition of residue symbols we then have

$$a \equiv t \,,\; a' \equiv t' \,,\; r \equiv 1 - y \,,\; r' \equiv 1 - y' \pmod p$$

The increase in the $p$-level of the length of the $p$-rim is given by

$$d = \left[\frac{a}{p}\right] - \left[\frac{a'}{p}\right] .$$

By the definition of $s$, we have

$$r - r' \equiv s \pmod p$$

and we first want to show that in fact

$$r - r' = s$$

holds. For this, we recall the general inequalities that consecutive columns in the Mullineux symbol have to satisfy:

(1) $\varepsilon \le r - r' < p + \varepsilon$

(2) $r - r' + \varepsilon \le a - a' < p + (r - r') + \varepsilon$

If $r - r' \neq s$, then by inequality (1) we must have:

$$\varepsilon = 1 \,,\; s = 0 \,,\; r - r' = p \,.$$

But by Theorem 3.1, this does not correspond to an admissible JS-extension. Hence $r - r' = s$ as claimed, and the inequality (2) now reads

(3) $s + \varepsilon \le a - a' < s + p + \varepsilon$.

We now consider the different types of extensions case-by-case.

**Case 1.** $\varepsilon' = \varepsilon = 1$.

Then $x + y \equiv x' + y' \equiv \alpha \pmod p$, and hence

$$s \equiv y' - y \equiv x - x' \pmod p \,.$$

This implies

$$a - a' \equiv t - t' \equiv x - x' + y' - y \equiv 2s \pmod p \,.$$

From (3) we then obtain immediately $a - a' = 2s$ and hence

$$d = \left[\frac{a}{p}\right] - \left[\frac{a'}{p}\right] = \left[\frac{t' + 2s}{p}\right] - \left[\frac{t'}{p}\right] = \left[\frac{t' + 2s}{p}\right]$$

as claimed in the Proposition.

**Case 2.** $\varepsilon' = 1$, $\varepsilon = 0$.
Here $a \equiv t = 0$, and (3) gives

$$s + 1 \leq a - a' < p + s + 1 \, .$$

Since $a - a' \equiv -t'$, we see that

$$a - a' = \begin{cases} p - t' & \text{if } s + t' \leq p - 1 \\ 2p - t' & \text{if } s + t' > p - 1 \end{cases}$$

and thus

$$d = \left[\frac{a}{p}\right] - \left[\frac{a'}{p}\right] = \left[\frac{t' + s + p}{p}\right] \, .$$

**Case 3.** $\varepsilon' = 0$.
Here we have $a' \equiv t' = 0$, and (3) gives

$$s \leq a - a' < p + s \, .$$

Since $a - a' \equiv t - t' = t$, we see that

$$a - a' = \begin{cases} t & \text{if } s \leq t \\ t + p & \text{if } s > t \end{cases}$$

and thus

$$d = \left[\frac{a}{p}\right] - \left[\frac{a'}{p}\right] = \left[\frac{p - 1 + s - t}{p}\right] \, .$$

$\diamond$

In the following Theorem we refine the description of JS-partitions via paths in a diagram (as given before) by also computing the $p$-weight of the partition along the path. We recall that for simplification we had combined some edges going into and out of the singular columns in the previous diagram; for the weight labeling we have to be more careful as these edges now have to be treated separately. The preceding proposition has already dealt with one possible contribution to the weight and thus explains one part of the label we now attach to an edge. The labeled diagram describes the weight changes along the graph for constructing JS-partitions of type $\alpha$; the labeling $\binom{d}{e}$ of an edge corresponds to an increase of $d$ in the weight resulting from an increase of the $p$-level of the $p$-rim length (this increases the weight from then on), respectively to an increase by $e$ (resulting from a possibly decreasing core) which adds on to the weight at this particular step only.

**Theorem 3.9** *To each edge* $\begin{smallmatrix} x' \\ y' \end{smallmatrix} \to \begin{smallmatrix} x \\ y \end{smallmatrix}$ *in the construction diagram of JS-partitions of type* $\alpha$ *we attach a label* $\binom{d}{e}$ *as follows.*

*Define* $t, t', s \in \{0, 1, \ldots, p-1\}$ *by* $t \equiv x - y + 1$, $t' \equiv x' - y' + 1$, $s \equiv y' - y$ *(mod p) as in the previous proposition, and let* $\varepsilon$ *and* $\varepsilon'$ *be as before. Furthermore, let* $\widetilde{\mu}$ *and* $\mu$ *be the p-cores associated to the paths ending at* $\begin{smallmatrix} x' \\ y' \end{smallmatrix}$ *and* $\begin{smallmatrix} x \\ y \end{smallmatrix}$*, respectively.*

*Then we put*

$$d = \begin{cases} [\frac{t'+2s}{p}] & \text{if } \varepsilon' = \varepsilon = 1 \\ [\frac{t'+s+p}{p}] & \text{if } \varepsilon' = 1, \varepsilon = 0 \\ [\frac{p-1+s-t}{p}] & \text{if } \varepsilon' = 0 \end{cases} \quad \text{and} \quad e = \begin{cases} 1 & \text{if } |\mu| < |\widetilde{\mu}| \\ 0 & \text{if } |\mu| \geq |\widetilde{\mu}| \end{cases} .$$

*Furthermore, we attach to the only admissible singular start in the diagram the weight label* $d_0 = 1$*, and to the admissible regular starts the weight label* $d_0 = 0$*.*

*Let* $\lambda$ *be a JS-partition constructed by starting at a vertex with weight label* $d_0$ *and then proceeding along a path in the diagram with edges labeled* $\binom{d_1}{e_1}, \ldots, \binom{d_k}{e_k}$*. Then the p-weight of* $\lambda$ *is*

$$w(\lambda) = \sum_{i=0}^{k} (k+1-i) d_i + \sum_{i=1}^{k} e_i \,.$$

**Proof.** Suppose that we build up a JS-partition $\lambda$ following a path with $k$ steps, i.e.

$$G_p(\lambda) = \begin{pmatrix} a_{k+1} & a_k & \cdots & a_1 \\ r_{k+1} & r_k & \cdots & r_1 \end{pmatrix} \quad \text{resp. } R_p(\lambda) = \begin{Bmatrix} x_1 & \cdots & x_k & x_{k+1} \\ y_1 & \cdots & y_k & y_{k+1} \end{Bmatrix}.$$

The start from the empty partition $\emptyset$ to $\binom{a_1}{r_1}$ (resp. $\left\{\begin{smallmatrix} x_1 \\ y_1 \end{smallmatrix}\right\}$) might be considered as the 0-th step, with label $d_0 = 1 - \varepsilon$ and $e_0 = 0$.

Let $m_{k+1}$ be the $p$-level of $a_{k+1}$, i.e.

$$a_{k+1} = m_{k+1} p + t_{k+1} \quad \text{with } t_{k+1} \in \{0, \ldots, p-1\}$$

The upper edge labels in the diagram given in the statement of the Theorem are chosen exactly according to the previous Proposition; thus by that Proposition the labels $d_0, \ldots, d_{k+1}$ describe the increase of the $p$-level at each step and hence

$$m_{k+1} = \sum_{i=0}^{k} d_i \,.$$

We now prove the claim on the weight by induction on the number $k$ of construction steps for $\lambda$. For $k = 0$ the assertion clearly holds by the definition of $d_0$.

Let $\widetilde{\lambda}$ be the partition obtained after the first $k-1$ steps, i.e.

$$G_p(\widetilde{\lambda}) = \begin{pmatrix} a_k & \cdots & a_1 \\ r_k & \cdots & r_1 \end{pmatrix} .$$

Then by induction we know that

$$w(\widetilde{\lambda}) = \sum_{i=0}^{k-1}(k-i)d_i + \sum_{i=1}^{k-1} e_i \,.$$

By Theorem 3.6 the $p$-core of a JS-partition constructed along a path in the graph given in Theorem 3.1 is determined by the $y$-entry of the column where the path ends. Now let $\widetilde{\mu}$ and $\mu$ be the $p$-cores associated to the paths ending at ${x_k \atop y_k}$ and ${x_{k+1} \atop y_{k+1}}$, respectively. By definition, the lower edge labels are chosen to be

$$e_k = \begin{cases} 1 & \text{if } |\mu| < |\widetilde{\mu}| \\ 0 & \text{if } |\mu| \geq |\widetilde{\mu}| \end{cases} \,.$$

As

$$|\mu| - |\widetilde{\mu}| \equiv |\lambda| - |\widetilde{\lambda}| = a_{k+1} \equiv t_{k+1} \pmod p$$

and $p$-cores associated with neighbouring vertices in the JS-construction diagram differ by less than $p$ boxes, we deduce

$$|\mu| - |\widetilde{\mu}| = t_{k+1} - e_k p \,.$$

By definition of the weight, this implies

$$a_{k+1} = |\lambda| - |\widetilde{\lambda}| = p(w(\lambda) - w(\widetilde{\lambda})) + |\mu| - |\widetilde{\mu}| = p(w(\lambda) - w(\widetilde{\lambda}) - e_k) + t_{k+1} \,.$$

Hence

$$w(\lambda) - w(\widetilde{\lambda}) - e_k = m_{k+1} = \sum_{i=0}^{k} d_i \,,$$

and using the induction hypothesis we conclude

$$w(\lambda) = \sum_{i=0}^{k}(k+1-i)d_i + \sum_{i=1}^{k} e_i$$

as was to be proved. $\diamond$

# 4 The Mullineux map and Mullineux fixed JS-partitions

The Mullineux map is an involutory bijection $\lambda \mapsto \lambda^M$ on the set of $p$-regular partitions. It is defined using the Mullineux symbols. If

$$G_p(\lambda) = \begin{pmatrix} a_1 & a_2 & \cdots & a_m \\ r_1 & r_2 & \cdots & r_m \end{pmatrix} ,$$

then

$$G_p(\lambda^M) = \begin{pmatrix} a_1 & a_2 & \cdots & a_m \\ s_1 & s_2 & \cdots & s_m \end{pmatrix}.$$

where $s_j = a_j - r_j + \varepsilon_j$, and

$$\varepsilon_j = \begin{cases} 0 & \text{if } p \mid a_j \\ 1 & \text{otherwise} \end{cases}.$$

An easy calculation shows that if

$$R_p(\lambda) = \left\{ \begin{matrix} x_1 & \cdots & x_m \\ y_1 & \cdots & y_m \end{matrix} \right\}$$

then the residue symbol of the Mullineux conjugate $\lambda^M$ is

$$R_p(\lambda^M) = \left\{ \begin{matrix} \delta_1 - y_1 & \cdots & \delta_m - y_m \\ \delta_1 - x_1 & \cdots & \delta_m - x_m \end{matrix} \right\}$$

where

$$\delta_j = \begin{cases} 1 & \text{if } x_j + 1 = y_j \\ 0 & \text{otherwise} \end{cases}.$$

Thus $\delta_j = 1 - \varepsilon_j$. From the above we see that the entries of the residue symbols for a Mullineux fixed point satisfy $x_j + y_j = \delta_j$ for all $j$.

We recall that a column $\substack{a_j \\ r_j}$ in $G_p(\lambda)$ is called *singular* if $p|a_j$, otherwise the column is called *regular*. The corresponding condition in the residue symbol is $x_j + 1 = y_j$. Thus for singular columns we have $\delta_j = 1$, and for regular columns $\delta_j = 0$.

Since the Mullineux map is the identity for $p = 2$ we will from now on always assume that $p > 2$.

We now give a description of those JS-partitions which are fixed under the Mullineux map.

**Theorem 4.1** *Let $\lambda$ be a $p$-regular partition, $p > 2$.*

*Then $\lambda$ is a JS-partition and a Mullineux fixed point if and only if $R_p(\lambda) = \left\{ \begin{matrix} x_1 & \cdots & x_k \\ y_1 & \cdots & y_k \end{matrix} \right\} \in \mathcal{R}_0^f$, which is defined to be the set of residue symbols constructed iteratively by the following procedure:*

$$\left\{ \begin{matrix} x_1 \\ y_1 \end{matrix} \right\} = \left\{ \begin{matrix} 0 \\ 0 \end{matrix} \right\}$$

*If the first $k-1$ columns of the residue symbol are already constructed then we have two possibilities for a regular extension, namely*

$$\begin{matrix} x_k \\ y_k \end{matrix} = \begin{matrix} 1 - y_{k-1} \\ y_{k-1} - 1 \end{matrix} \quad \text{or} \quad \begin{matrix} y_{k-1} - 1 \\ 1 - y_{k-1} \end{matrix}$$

*and if $y_{k-1} = 1$, then we have furthermore one possibility for a singular extension, namely*

$$\begin{matrix} x_k \\ y_k \end{matrix} = \begin{matrix} 0 \\ 1 \end{matrix} .$$

*In particular, Mullineux fixed JS-partitions are always of type 0.*

**Proof.** This follows easily from the construction rules for the residue symbols of JS-partitions together with the added restrictions coming from the Mullineux fixed point condition. $\diamond$

Again, the construction rule given in the theorem can be described alternatively by a diagram in the following way. The JS-partitions which are Mullineux fixed points are obtained by adding on columns to the residue symbol on a walk starting at $\begin{matrix} 0 \\ 0 \end{matrix}$ in the diagram $(*)$ below (again, we omit the brackets for simplification, and the symbol $\subset$ in the diagram means that we have a loop at the corresponding node of the graph). The edges in the diagram are labelled $\binom{d}{e}$ according to Theorem 3.9; labels $\binom{0}{0}$ are omitted.

$(*)$

$$\begin{array}{ccccccccccc}
 & \nearrow & \boxed{\begin{matrix} 0 \\ 0 \end{matrix}} & \rightarrow & \begin{matrix} 1 \\ -1 \end{matrix} & \rightarrow & \begin{matrix} 2 \\ -2 \end{matrix} & \rightarrow & \cdots & \rightarrow & \begin{matrix} \frac{p-3}{2} \\ -\frac{p-3}{2} \end{matrix} \\
\begin{matrix} \bigcap \\ 0 \\ 1 \end{matrix} & & \binom{1}{1} \downarrow\uparrow \binom{1}{0} & & \binom{1}{1} \downarrow\uparrow \binom{1}{0} & & \binom{1}{1} \downarrow\uparrow \binom{1}{0} & & \cdots & & \binom{1}{1} \downarrow\uparrow \binom{1}{0} \\
 & \underset{\binom{1}{0}}{\nwarrow} & \begin{matrix} -1 \\ 1 \end{matrix} & \overset{\binom{0}{1}}{\leftarrow} & \begin{matrix} -2 \\ 2 \end{matrix} & \overset{\binom{0}{1}}{\leftarrow} & \begin{matrix} -3 \\ 3 \end{matrix} & \overset{\binom{0}{1}}{\leftarrow} & \cdots & \overset{\binom{0}{1}}{\leftarrow} & \begin{matrix} -\frac{p-1}{2} \\ \frac{p-1}{2} \end{matrix}
\end{array}$$

This is clearly a subgraph of the previous graph for JS-partitions of type 0, where all singular columns except $\begin{matrix} 0 \\ 1 \end{matrix}$ are deleted since they are not fixed under the Mullineux map.

As we have seen in §3 JS partitions of type $\alpha$ have empty or rectangular $p$-cores with corner node of residue $\alpha$; thus Mullineux fixed JS-partitions have empty or square $p$-cores. More precisely, the previous Theorem tells us how to compute the $p$-core from the final column in the residue symbol, i.e. how it is determined by the end of the path in the graph.

**Theorem 4.2** *The $p$-core of a Mullineux fixed JS-partition with residue symbol $R_p(\lambda) =$*

$\left\{\begin{array}{ccc} x_1 & \cdots & x_m \\ y_1 & \cdots & y_m \end{array}\right\}$ *is given as follows.*

$$\lambda_{(p)} = \begin{cases} (l^l) \quad \text{if} \quad \begin{array}{c} x_m \\ y_m \end{array} \in \left\{ \begin{array}{c} l-1 \\ 1-l \end{array}, \begin{array}{c} \alpha-l-1 \\ l+1 \end{array}, \begin{array}{c} l \\ 1+l \end{array}, \begin{array}{c} -l \\ 1-l \end{array} \right\} \\ \qquad \text{and } 1 \leq l \leq \frac{p-1}{2} \\ \emptyset \quad \text{if} \quad \begin{array}{c} x_m \\ y_m \end{array} \in \left\{ \begin{array}{c} -1 \\ 1 \end{array}, \begin{array}{c} 0 \\ 1 \end{array} \right\} \end{cases}$$

For a better illustration, we have put the $p$-cores at the corresponding end points of the walks in the diagram:

$$\begin{array}{ccccccccccc} & \nearrow & \boxed{(1)} & \to & (2^2) & \to & (3^3) & \to & \cdots & \to & (\frac{p-1}{2}^{\frac{p-1}{2}}) \\ \bigcap & & & & & & & & & & \\ \emptyset & & \downarrow\uparrow & & \downarrow\uparrow & & \downarrow\uparrow & & \cdots & & \downarrow\uparrow \\ & \nwarrow & \emptyset & \leftarrow & (1) & \leftarrow & (2^2) & \leftarrow & \cdots & \leftarrow & (\frac{p-3}{2}^{\frac{p-3}{2}}) \end{array}$$

We have seen before that the existence of JS-partitions for a given rectangular $p$-core and weight was easy to answer in the affirmative.

It is more complicated to answer the question about the existence of Mullineux fixed JS-partitions for a given square $p$-core and even weight (note that the weight of any Mullineux fixed $p$-regular partition is even [2, Prop. 3.4]).

**Theorem 4.3** *Let $w$ be an even weight, $\mu$ a square $p$-core or empty. Then there always exists a Mullineux fixed JS-partition of weight $w$ and with $p$-core $\mu$, except in the case $w = 2$ and $\mu = (\frac{p-1}{2}^{\frac{p-1}{2}})$.*

**Proof.** We consider the weight changes along the paths in the diagram $(*)$ (given after Theorem 4.1) as described by Theorem 3.9.

For each square $p$-core $\mu$ and even weight $w$ we have to describe a path in the diagram $(*)$ ending in a vertex with the prescribed core associated. We abbreviate the vertices by writing only the top residue in the way it is given in $(*)$, i.e. instead of $\begin{array}{c} a \\ -a \end{array}$ we only write $a$, and we write $0'$ for the singular column $\begin{smallmatrix} 0 \\ 1 \end{smallmatrix}$.

**Case $\mu = \emptyset$.**

For $w = 0$, we just take the empty partition.

Now let $w$ be a positive even number. The path in $(*)$ then has to end at one of the vertices $-1$ or $0'$.

For $w = 2$, we take the path $0 \to -1$.

For $w \geq 4$, we choose the path $0 \to -1 \to (0')^m$, meaning that the path ends on $m-1 = \frac{w-4}{2}$ loops at $0'$.

**Case** $\mu = (j^j)$, $j \in \{1, \ldots, \frac{p-3}{2}\}$**.**
In this case the path has to end at one of the vertices $j-1$ or $-(j+1)$.
For $w = 0$, we take the path $0 \to 1 \to \cdots \to j-2 \to j-1$.
For $w = 2$, we take the path $0 \to 1 \to \cdots \to j \to -(j+1)$.
For $i \in \{0, \ldots, j-1\}$ the path

$$0 \to 1 \to \cdots \to i \to -(i+1) \to i \to i+1 \to \cdots \to j-1$$

constructs a partition as required of weight $2(j+1-i)$. This takes care of the weights $w \in \{4, 6, \ldots, 2j+2\}$.
For $w \geq 2j+4$, we choose the path

$$0 \to -1 \to (0')^m \to 0 \to 1 \to \cdots \to j-1$$

with $m-1 = \frac{w-2j-4}{2}$ loops at $0'$, along which a JS-partition of weight $w$ as required is constructed.

**Case** $\mu = (\frac{p-1}{2}^{\frac{p-1}{2}})$**.**
Here the paths above may also be chosen, except in the case $w = 2$, where the path given above would not end at a vertex associated with the correct core $(\frac{p-1}{2}^{\frac{p-1}{2}})$. In this case, no path to the only vertex $\frac{p-3}{2}$ with the correct core yields $w = 2$, because the first edge giving a positive weight is a vertical one down, and to get back to the vertex $\frac{p-3}{2}$ (which is at the end of the top row) an additional weight is added, giving a weight $> 2$. $\diamond$

**Acknowledgements**
The authors gratefully acknowledge the support of the Danish Natural Science Foundation. The final version of this article was written during a very pleasant stay at the Mathematical Sciences Research Institute at Berkeley; the research at MSRI was supported in part by NSF grant # DMS 9022140 as well as by the Deutsche Forschungsgemeinschaft (grant Be 923/6-1 of the first author). We would like to thank these institutions for their support.

Christine Bessenrodt, Fakultät für Mathematik, Otto-von-Guericke-Universität Magdeburg, 39016 Magdeburg, Germany
*Email address:* bessen@mathematik.uni-magdeburg.de

Jørn B. Olsson, Matematisk Institut, Københavns Universitet, Universitetsparken 5, 2100 Copenhagen Ø, Denmark
*Email address:* olsson@math.ktu.dk